\documentclass{amsart}
\usepackage[T1]{fontenc}
\usepackage[utf8]{inputenc}
\usepackage{amssymb}
\usepackage{amsmath}
\usepackage{graphicx}
\usepackage{rotating}
\def\R{I\kern -0,37 em R}
\def\P{I\kern -0,37 em P}
\def\Z{I\kern -0,37 em Z}
\begin{document}
\title{Automorphisms of Flag Systems}

\author{A. Kumpera}
\address{State University of Campinas, Campinas SP, Brazil. Tel.: +55-44-32253037}
\email{xyz@math.university.edu}

\maketitle

\begin{abstract}
We first discuss the problems in the theory of ordinary differential equations that gave rise to the concept of a flag system and illustrate these with the Cartan criterion for Monge equations (1st order) as well as the Cartan statement concerning the local equivalence of Monge-Ampère type equations (2nd order). Next, we describe a prolongation functor operating on the infinitesimal symmetries (automorphisms) of the Darboux flag and extending these, isomorphically, to all the symmetries of any other flag. Hence, flag systems cannot be distinguished by their symmetry algebras and the local classification of these objects is approached by considering higher order isotropies of these algebras as well as the groupoids of $k-th$ order formal equivalences since the differential equations defining the latter provide precious information for the application of flag systems to differential equations (\textit{e.g.}, Cartan's criterion for non-linear Monge equations). In examining the behaviour of the isotropy algebras, that can either diminish or remain the same, when passing from a derived system $S_\nu$ to the previous system $S_{\nu-1}$, we obtain a full set of numerical invariants for the elementary flag systems that moreover specify the local models.
\end{abstract}

\section{Introduction}
Special Pfaffian systems, called \textit{flags systems} or simply \textit{flags}, are very interesting and important non-integrable systems. They are characterised by the property that their associated sequence of \textit{derived} systems is the longest possible (its length being equal to the rank of the Pfaffian system) hence the non-integrability property manifests itself "in slow steps" (to be clarified later). This being so, we have at our disposal several techniques and methods that make such systems among the most relevant in Geometric Control Theory namely, in the study of accessibility and controllability. However, perhaps the most remarkable results are, as evidenced by Élie Cartan, those related to the integrability of systems of ordinary and partial differential equations. To each such system is associated a contact Pfaffian system and the study of these contact systems as well as those associated to their prolongations is the magical touch given by Cartan to the integration problem. We would also like to mention that these flag systems might be of considerable interest in cosmology since the behavior of the uncountably many flags with respect to the "expansion" of the Darboux automorphisms action, tending towards a free action, has rather striking similarities with the expansion of the universe, including the behaviour of its mega-singularities (black holes). Finally, it seems worthwhile to mention that flag systems are \textit{rigid} structures in the sense that they cannot be continuously deformed along equivalent systems and, in fact, each flag has its own specific \textit{model}.  

Flag systems go back to Gaston Darboux and Friedrich Engel \cite{Pfaff1814,Darboux1882,Engel1889,Monge1889}. The well known Darboux Theorem \cite{Kumpera1970} provides the local description of \textit{length one flags} and the Engel Theorem the description of \textit{length two flags} \cite{Cartan1899,Cartan1901,Engel1889,Kumpera1982}. Élie Cartan described the local model for a \textit{homogeneous or transitive} flag and used this object to study the solutions of under-determined systems of ordinary differential equations. He gave a criterion under which the general solution of such a system can be given an explicit form evidencing its dependence on an arbitrary function of one variable and on all its derivatives up to a finite order $\mu$ \cite{Cartan1914}, see also \cite{Goursat1905,Kumpera1999,Zervos1913}.

Édouard Goursat also gave attention to flag systems and, in the preface of his book \cite{Goursat1922}, he candidly writes

\vspace{3 mm}

\noindent 
\textit{Enfin, dans les trois derniers Chapitres, j'expose quelquesuns des progrès les plus récents acquis à la science, relatifs aux systèmes de Pfaff. Les plus importants de ces progrès sont dus à M. Cartan, dont on trouvera le nom presque à chaque page de ces trois chapitres.}

\vspace{3 mm}

In chapter VII, p.328, \textit{l.}26-29, Goursat describes the canonical form (\textit{forme canonique}) to which special systems (his terminology for flag systems) can be reduced. Unfortunately, his statement is inexact since he only exhibits the Cartan homogeneous model. In much the same way, Engel and von Weber \cite{Weber1898} also ignore non-homogeneity. To our knowledge, Cartan never mentioned explicitly any non-homogeneous model but he was certainly aware of them as we shall see later. 

\vspace{2 mm}

As for the presentation, we try to be as concise as possible by providing all the required references. Nevertheless, we are forced to recall some well known concepts, definitions and results basically due to the fact that we are inclined, along with Élie Cartan and Édouard Goursat, to work in the realm of "covariance" whereas most authors (perhaps all) prefer the "contravariant" setting. Concerning the examples and the (extensive) calculations, they are to be considered not only as illustrations but, in fact, play a fundamental role in the development 
of the theory since they are the \textit{foundations} as well as the \textit{pillars} for the flag systems. 
All flags have a common foundation provided by the Darboux flag and a common first pillar namely the Engel flag. 
As for the second pillar, we are faced with two options namely, the Cartan homogeneous (transitive) flag and the 
rather bizarre exceptional (non-transitive) flag for which, amazingly, the Cartan criterion for Monge equations 
(only provable in the homogeneous situation) does not break down in spite of a co-dimension one singular locus 
\cite{Zervos1932,Kumpera1987}. A closer inspection reveals that Cartan achieved this remarkable feat by examining the structure of the algebra of infinitesimal automorphisms of the exceptional flag in five dimensional space 
and its isotropy up to second order. In a later section, we also examine other (exponentially proliferating) 
pillars. As for the extensive calculations, both Lie and Cartan frequently indulged into incredible calculations 
since they believed, presumably, that this was the first step towards understanding Heaven right here from earth. 
Examining our calculations within the Hamiltonian formalism, especially those concerning the isotropy algebras, it is worthwhile to observe the amazing (and meaningful) indexation assumed by the partial derivatives.

\vspace{2 mm}

\noindent
\textbf{ADDENDUM}. Very recently, the author found a mirror image of the result stated in the theorem below in Mormul's \textit{Hokkaido Math. J.} \textbf{34}(2005)1-35 on page 17, \textit{l.} -11...-5, though in a very different context. His is Wagnerian in seriousness (\textit{Nibelungen}) whereas ours is Mozartian in joy (\textit{the flute concertos}). Moreover, it seems to the author that Mormul's \textit{mériédrique} argument, in the sense of Élie Cartan, is not altogether appropriate. Anyway, all this is completely irrelevant since Élie Cartan provided, more than a century ago, a full account of this result as well as many other results that both authors most probably ignore.
 
\section{Flag Systems}
A Pfaffian system on a manifold \textit{M} (of \textit{dim n}) will be envisaged as a locally trivial vector sub-bundle  $S\subset T^*M=T^*$.
For simplicity, we assume that all the data is of class $C^\infty$ and that \textit{M} is connected and second countable. The (constant) rank of \textit{S} is denoted by \textit{r} and we shall say, occasionally, that a linear differential form $\omega$ belongs to \textit{S} whenever $\omega\in\Gamma(S)$, the later denoting the pre-sheaf of all local sections of \textit{S}. We recall that \textit{S} is integrable when, for all $\omega\in\Gamma$(\textit{S}),
\begin{equation}
d\omega\equiv 0~mod~S\hspace{7 mm}i.e.,\hspace{7 mm}d\omega=\sum\sigma^i\wedge\omega^i,\hspace{15 mm}\omega^i\in\Gamma(S)
\end{equation}

\noindent
The Martinet \textit{structure tensor} \cite{Kumpera1982,Martinet1974}, is the vector bundle morphism
\begin{equation}
\delta :S\longrightarrow\wedge^2(T^*/S)
\end{equation}

\noindent 
defined on local sections by $\delta(\omega)=d\omega~mod~S$. Since 
\begin{equation*}
d(f\omega)=fd\omega+df\wedge\omega\equiv fd\omega~mod~S
\end{equation*}

\noindent
the above defined pre-sheaf morphism is linear over the functions of \textit{M} and therefore induces the desired vector bundle morphism. We now assume that the rank of $\delta$ is constant and define the \textit{derived system} $S_1=ker~\delta$. Then $S_1$ is also a Pfaffian system and $S_1\subset S$. The integrability condition (1) now reads $S_1=S$ and we can say, intuitively, that a non-integrable system \textit{S} deviates least from integrability when $rank~S_1=r-1$. In fact, it means that \textit{most} local sections of \textit{S} do satisfy that integrability condition. Furthermore, we can transcribe the definition of $S_1$ in a more practical way, providing an operational procedure for its calculation. 
For this, we simply observe that the two equivalent conditions stated in (1) apply as well, independently of any concern about integrability, to an arbitrary local section $\omega$ of \textit{S} and where the forms $\omega^i$ are also any local sections of \textit{S}. Translated in a more academic language, we are just saying that $S_1$ is the Pfaffian system generated by all the local forms $\omega\in\Gamma(S)$ such that $d\omega$ belongs to the pre-sheaf of ideals, in $\Gamma(\wedge T^*)$, generated by $\Gamma(S)$ and where $\Gamma(~)$ denotes the set of local sections.
We can now define inductively the $(\nu+1)-st$ derived system $S_{\nu+1}=(S_\nu)_1$ by considering at each step the Martinet tensor (2) adapted to the system $S_\nu$ and assuming of course that the rank of this tensor is constant. This iterated construction yields a decreasing sequence of Pfaffian systems that forcefully becomes stationary, namely:

\begin{equation}
S=S_0\supset S_1\supset S_2\supset\cdots\supset S_\ell=S_{\ell+1}
\end{equation}

\vspace{2 mm}

\noindent
The integer $\ell$ is called the \textit{length} of the Paffian system \textit{S}. The last term $S_\ell$ in the above sequence is integrable and we shall refer to it as the \textit{terminating} system of the flag. Length zero characterises the integrability of \textit{S}. A Pfaffian system for which all the successive structure tensors are of constant rank is called \textit{totally regular} \cite{Libermann1978}. We say that \textit{S} is a \textit{flag system} (Système de Pfaff en drapeau, \cite{Kumpera1982} when $rank~S_\nu=r-\nu,~\nu\leq\ell$. The basic theory as well as an initial attempt towards a local classification of flag systems can be found in the previous citation.

Flag systems and more generally non-integrable Pfaffian systems are extremely useful in the analysis of systems of ordinary and partial differential equations. The pioneering work is entirely due to Sophus Lie culminating in his \textit{Opera Magna} named \textit{Verwertung des Gruppenbegriffes für Differentialgleichungen} \cite{Lie1895}. Élie Cartan also gave extraordinary contributions to this subject and no doubt went fathoms deeper. Among the several papers he wrote, perhaps the deepest and above all the most inspiring are \cite{Cartan1910} and \cite{Cartan1914} (see also the forewords in \cite{Kumpera1999}, \cite{Kumpera2000} and \cite{Kumpera2002}. Knowing the structure and, whenever possible, a local model for the Pfaffian system \textit{S} associated to the set of equations, we shall be able \cite{Kumpera1972} to determine the algebra of all infinitesimal automorphisms (dwelling in the base space) of the given $k-th$ order equations, lift this algebra to a sub-algebra of infinitesimal $k-th$ order contact transformations and subsequently prolong the later to an order where the equations become formally integrable (\textit{e.g}., involutive).The prolongued algebra is then equal to the set of all infinitesimal contact automorphisms of the prolongued equations and \textit{a fortiori} leaves invariant the characteristic system of the prolongation of \textit{S}. The knowledge of such an algebra entails the possibility of integrating, via Lie and Cartan's methods, the (integrable) characteristic system and thereafter the initially given equations \cite{Kumpera1999,Kumpera2002}.

Élie Cartan studies, in \cite{Cartan1910}, second order \textit{involutive} systems of two equations, namely the Monge equations \cite{Kumpera1999}, and single second order equations with integrable double Monge characteristics. The flag systems of length three appear in connection with certain types of \textit{linear} and involutive systems of two equations \cite{Cartan1910}, p.126, II). These flags are always homogeneous and it turns out that the linear systems are invariant under an algebra of second order infinitesimal contact transformations obtained by lifting and subsequently prolonging the algebra of all infinitesimal transformations leaving invariant the (Pfaffian) equation $dx^2+x^3dx^1=0$ or, in other words, the Darboux system in the space of three variables. In the case of second order equations with integrable double Monge characteristics, we can single out the Monge-Ampère equation \cite{Goursat1891,Goursat1898,Goursat1905}

\begin{equation*}
(r+A)(t+C)-(s+B)^2=0
\end{equation*}

\noindent
where \textit{r,s} and \textit{t} denote the second order partial derivatives (\textit{s} being the mixed derivative). The associated Pfaffian system is a flag system of rank and length two and its \textit{covariant} system, of rank three, is therefore integrable \cite{Kumpera1999,Kumpera2002}. Cartan shows that all the equations with such an associated Pfaffian system are locally equivalent, the above Monge-Ampère equation being the standard model, and that the algebra of all infinitesimal second order contact automorphisms of these equations is (locally) isomorphic to the algebra of all vector fields on $\bf{R}^3$. We infer that the above associated Pfaffian system is the Engel flag. Another noteworthy example is given by the second order Goursat equations \cite{Goursat1898,Kumpera1999,Kumpera2000} whose algebra of all infinitesimal second order contact automorphisms is isomorphic to the exceptional complex simple Lie algebra $g_2$ with real form $g_{2(2)}$.

\section{Pseudo-normal forms}
We now recall the "approximate" pseudo-normal form of writing down, locally, a flag system in the sense that it will eventually lead to the discovery of true models and, further, enable to investigate the nature of infinitesimal automorphisms and isotropies. Nevertheless, the complicated nature of the write up of such pseudo-models will also lead us to esoteric descriptions of the aforementioned objects. It should be noted that 
the same flag can have several pseudo-normal forms. However, if two flags have the same pseudo-normal form, they are of course equivalent.

\newtheorem{pseudo}[LemmaCounter]{Lemma}
\begin{pseudo}
Let S be a flag system of length $\ell$~=~n defined on an $(n+2)$-dimensional manifold M, $n\geq 1$, and terminating in $0$. Then, in a neighborhood of each point $p\in M$, the system S can be expressed by the following set of generators, where the local coordinates $(x^1,x^2,~\cdots~,x^{n+2})$ are centered at the origin $0\in\bf{R^{n+2}:}$

\begin{equation*}
\hspace{5 mm}\omega^1=dx^{i_1}+~x^3dx^{j_1},\hspace{10 mm}(i_1,j_1)=(2,1)\hspace{40 mm}
\end{equation*}
\begin{equation*}
\hspace{5 mm}\omega^2=dx^{i_2}+~x^4dx^{j_2},\hspace{10 mm}(i_2,j_2)=(3,1)\hspace{40 mm}
\end{equation*}
\begin{equation*}
\hspace{3 mm}\omega^3=dx^{i_3}+~x^5dx^{j_3},\hspace{10 mm}(i_3,j_3)\in\{(4,j_2),(j_2,4)\}\hspace{23 mm}
\end{equation*}
\begin{equation*}
\hspace{3 mm}\omega^4=dx^{i_4}+X^6dx^{j_4},\hspace{10 mm}(i_4,j_4)\in\{(5,j_3),(j_3,5)\}\hspace{23 mm}
\end{equation*}
\begin{equation*}
\hspace{3 mm}\omega^5=dx^{i_5}+X^7dx^{j_5},\hspace{10 mm}(i_5,j_5)\in\{(6,j_4),(j_4,6)\}\hspace{23 mm}
\end{equation*}
\begin{equation*}
................................
\end{equation*}
\begin{equation*}
\omega^\nu=dx^{i_\nu}+X^{\nu+2}dx^{j_\nu},\hspace{6 mm}(i_\nu,j_\nu)\in\{(\nu+1,j_{\nu-1}),(j_{\nu-1},\nu+1)\}
\end{equation*}
\begin{equation*}
................................
\end{equation*}
\begin{equation*}
~\omega^n=dx^{i_n}+X^{n+2}dx^{j_n},\hspace{6 mm}(i_n,j_n)\in\{(n+1,j_{n-1}),(j_{n-1},n+1)\}
\end{equation*}

\vspace{6 mm}
\noindent
where, for $~6\leq\nu+2\leq n+2~,~X^{\nu+2}=x^{\nu+2}$ when $(i_\nu,j_\nu)=(j_{\nu-1},\nu+1)$ and $X^{\nu+2}=x^{\nu+2}+c^{
\nu+2}$ when $(i_\nu,j_\nu)=(\nu+1,j_{\nu-1})$. The coefficients $c^6,c^7,\cdots,c^{n+2}$ are real constants.

\noindent
The iterated derived systems of S, in the neighborhood of p, are obtained by removing successively the bottommost generators.

\noindent
If $j_3=j_4=~\cdots~=j_\nu=1~(i.e.,~i_3=4,~i_4=5,~\cdots~,~i_\nu=\nu+1)$ for a certain value $4\leq \nu\leq n$, the constants $c^6,~c^7,~\cdots~,~c^{\nu+2}$ can be chosen null, the expressions of the forms $\omega^1,~\omega^2,~\cdots~,\omega^{\nu+1},~\cdots~,~\omega^n$ remaining unchanged. When $\nu=n$, we obtain the Cartan transitive model.
\end{pseudo}

\vspace{4 mm}
\noindent
Proofs can be found in \cite{Kumpera1982} and \cite{Mormul2000}. We observe however that our pseudo-normal form slightly differs from the usual contravariant presentations, not only in its setting but mainly due to the single indexation in contra-position to the double indexations found elsewhere which makes the "dualisation functor" behave somewhat unfaithfully. Just for reference and also for obvious reasons we shall call them the 
\textit{Cartan-Fibonacci forms}.

\noindent
A pseudo-normal form is called \textit{elementary} when all the constants $c^{\nu+2}$ vanish. It should be observed that this terminology is only meaningful in what touches the local representation itself and is not an invariant notion for the generated flag since such a system can have simultaneously \textit{elementary and non-elementary} pseudo-normal forms. The flag system \textit{S} defined (generated) on $~\bf{R}^{\ell+2}~$ by a given pseudo-normal form is called a \textit{pseudo-model}.

\vspace{2 mm}

\newtheorem{model}[LemmaCounter]{Lemma}
\begin{model}
Any pseudo-model is equivalent, via a diffeomorphism $\varphi$ of $~\bf{R}^{\ell+2}$, to a pseudo-model with an
elementary pseudo-normal form.
\end{model}

\noindent
It suffices, of course, to take as new coordinates $\overline{x}^{\nu+2}=x^{\nu+2}+c^{\nu+2}$. However, the resulting diffeomorphism $\varphi:\bf{R}^{\ell+2}\longrightarrow\bf{R}^{\ell+2}$ will not preserve, in general, the origin.

\vspace{2 mm}

\newtheorem{merde}[LemmaCounter]{Lemma}
\begin{merde}
For any fixed length $\ell$, two pseudo-models with different elementary forms are not locally equivalent in any neighborhood of the origin and their small growth vectors \cite{Mormul2000} as well as the ranks of their reduced tensors and reduced  nilpotent algebras \cite{Kumpera1982} are distinct.
\end{merde}

\section{Models in low dimensions}

All the local models discussed in this section are defined in $\bf{R^n}$, for some integer $\bf{n}$, and should be viewed as germs of flag systems at the origin. We shall now assume, as usual, that:

\vspace{3 mm}

\hspace{1 mm}(i) the \textit{class} of \textit{S} is equal to \textit{dim M},

\vspace{2 mm}

(ii) $S_\ell=0$\hspace{2 mm} i.e., the terminating system is null.

\vspace{2 mm}

\noindent
The first assumption is quite natural since there is no sense in dragging along useless variables (cf. \cite{Kumpera1982}, sect.3) for we can always factor (locally) the given system modulo its Cauchy characteristics. As for the second, it is merely a non-restrictive convenience \cite[Thm.4.2]{Kumpera1982}. Thereafter, we can also derive the following properties \cite[Corol.3.1]{Kumpera1982}:

\vspace{3 mm}

(iii) every flag system of length 1 is an odd class system of rank 1,

\vspace{2 mm}

(iv) the class of $S_\nu$ is equal to $n-\nu$ where $r=n-2$,

\vspace{2 mm}

\hspace{1 mm}(v) every flag system of length $\ell\geq 2$ has co-rank 2 hence its annihilator $\Sigma=S^\perp$ in \textit{TM} is a field of 2-planes. In particular, $dim~M=rank~S+2$ and $\ell=rank~S=length~S$.

\vspace{3 mm}

\noindent
$\bf{Flags~of~length~1 - Darboux}$

\vspace{2 mm}

\noindent
These are the Pfaffian systems of rank 1 and maximum odd class $2p+1=\ell+2$ hence, by the Darboux Theorem, have the local generator:
\begin{equation}
\omega=dx^0+\sum~x^idx^{p+i},\hspace{5 mm}1\leq i\leq p.
\end{equation}

\vspace{2 mm}

\noindent
We observe that non-integrable systems of rank 1 and constant class greater than 3 are also flag systems of length 1 but do not satisfy the assumption (i).
.
\vspace{2 mm}

\noindent
$\bf{Flags~of~length~2 - Engel}$

\vspace{2 mm}

\noindent
All such systems are locally equivalent to the Engel model in $\bf{R^4}$ \cite{Engel1889,Kumpera1982,Weber1898}.
\begin{equation*}
\omega^1=dx^2+x^3dx^1
\end{equation*}
\begin{equation*}
\omega^2=dx^3+x^4dx^1
\end{equation*}

\vspace{2 mm}

\noindent
The derived system is the Darboux flag generated by the form $\omega^1$, considered as a model in $\bf{R^3}$ (after reduction modulo the characteristics).

\vspace{2 mm}

\noindent
$\bf{Flags~of~length~3 - Cartan:~homogeneous~and~exceptional}$

\vspace{2 mm}

\noindent
There exist 2 local models, in $\bf{R^5}$, namely the homogeneous (transitive) Cartan model and a non-homogeneous (intransitive) exceptional model \cite{Giaro1978}:
\begin{equation*}
\omega^1=dx^2+x^3dx^1\hspace{15 mm}\varpi^1=dx^2+x^3dx^1
\end{equation*}
\begin{equation*}
\omega^2=dx^3+x^4dx^1\hspace{15 mm}\varpi^2=dx^3+x^4dx^1
\end{equation*}
\begin{equation*}
\omega^3=dx^4+x^5dx^1\hspace{15 mm}\varpi^3=dx^1+x^5dx^4
\end{equation*}

\vspace{2 mm}

\noindent
The derived systems of both models are equal to the Engel flag and the second model is non-homogeneous along the hyperplane $x^5=0$ ,  homogeneous elsewhere.

\vspace{2 mm}

\noindent
$\bf{Flags~of~length~4 - Cartan:~homogeneous~and~exceptionals}$

\vspace{2 mm}

\noindent
There exist $5=2^2+1$ non-equivalent models of length four defined in $\bf{R^6}$ namely:

\begin{equation*}
dx^2+x^3dx^1\hspace{18 mm}dx^2+x^3dx^1
\end{equation*}
\begin{equation*}
dx^3+x^4dx^1\hspace{18 mm}dx^3+x^4dx^1
\end{equation*}
\begin{equation*}
dx^4+x^5dx^1\hspace{18 mm}dx^4+x^5dx^1
\end{equation*}
\begin{equation*}
dx^5+x^6dx^1\hspace{18 mm}dx^1+x^6dx^5
\end{equation*}

\vspace{4 mm}

\begin{equation*}
dx^2+x^3dx^1\hspace{14 mm}dx^2+x^3dx^1\hspace{14 mm}dx^2+x^3dx^1
\end{equation*}
\begin{equation*}
dx^3+x^4dx^1\hspace{14 mm}dx^3+x^4dx^1\hspace{14 mm}dx^3+x^4dx^1
\end{equation*}
\begin{equation*}
dx^4+x^5dx^1\hspace{14 mm}dx^1+x^5dx^4\hspace{14 mm}dx^1+x^5dx^4
\end{equation*}
\begin{equation*}
dx^5+x^6dx^4\hspace{14 mm}dx^5+(x^6+1)dx^4\hspace{5 mm}dx^4+x^6dx^5
\end{equation*}

\vspace{4 mm}

\noindent
and $13=2^3+5$ models of length five in $\bf{R^7}$. Among the 14 models listed in \cite{Kumpera1982}, the 9th and 10th are equivalent \cite{Gaspar1985}. Further models can be found in \cite{Cheaito1999} and \cite{Mormul2000}. The general homogeneous flag of length $\ell$ has the Cartan local left model and the exceptional flag, with one frontal inversion, the right model that will be referred to as \textit{primary} \cite{Kumpera1982}.

\begin{equation*}
dx^2+x^3dx^1\hspace{18 mm}dx^2+x^3dx^1
\end{equation*}
\begin{equation*}
dx^3+x^4dx^1\hspace{18 mm}dx^3+x^4dx^1
\end{equation*}
\begin{equation*}
..................\hspace{19 mm}..................
\end{equation*}
\begin{equation*}
\hspace{4 mm}dx^\ell+x^{\ell+1}dx^1\hspace{15 mm}dx^\ell+x^{\ell+1}dx^1
\end{equation*}
\begin{equation*}
\hspace{4 mm}dx^{\ell+1}+x^{\ell+2}dx^1\hspace{15 mm}dx^1+x^{\ell+2}dx^{\ell+1}
\end{equation*}

\section{Infinitesimal automorphisms and isotropies}

Several numerical invariants and even letters have been considered for the purpose of classifying flag systems and exhibiting local models. Our purpose here is not exactly to search for additional or alternative invariants or to inquire for further models but to study the structure of the symmetry (invariance) groups and algebras associated to flag systems that, according to Lie and Cartan, are the magic key towards their understanding. All our considerations are, though sometimes veiled, of local nature and we now introduce "multi-fibrations" so as to improve our veiling methods. As usual, a \textit{fibration} is a surjective map $P\longrightarrow M$ of maximum rank and we define a \textit{multi-fibration} as being a sequence

\begin{equation}
P_\ell\longrightarrow P_{\ell-1}\longrightarrow\cdots\longrightarrow P_0 ~,
\end{equation}

\vspace{1 mm}

\noindent
where all the projections $\rho_{\beta,\alpha}:P_\alpha\longrightarrow P_\beta$ , $\alpha\geq\beta$, are fibrations. Given a Paffian system \textit{S} defined on a manifold \textit{M} and having regular successive derived systems, with regular characteristics, up to the final order $\ell$, we can take a suitable neighborhood $P_\ell$ of any point in \textit{M} and, restricting \textit{S} to this neighborhood as well as factoring the successive derived systems $S_\nu$, modulo their respective characteristic foliations, we can at least locally re-accommodate the given system \textit{S} as follows: $S=S_0$ still remains on $P_\ell$ (or on a quotient, modulo its characteristics) whereas each successive system $S_\nu$ finds a more adequate place on $P_{\ell-\nu}$ where it acquires maximum Cartan \textit{class}.

Without seeking any sort of pedantic "over precision", it seems however interesting, at least for the moment, to adopt some additional notations that later will be abandoned. We start by listing a few properties inherent to "naturality" and adopt the following abbreviations: \textbf{i.a.} stands for the \textit{infinitesimal automorphisms (symmetries)} of \textit{S}, $\mathcal{L}(S)$ for the Lie algebra of all these symmetries and finally $\chi(S)$ for the Cartan \textit{characteristic system} of \textit{S}.

\vspace{3 mm}

\hspace{2 mm}(i) All the systems $S_\nu$ are defined on the same space $P_\ell$ where $S=S_0$ is defined.

\vspace{2 mm}

\hspace{1 mm}(ii) Any \textbf{i.a.} of $S_\nu$ is also an \textbf{i.a.} of $S_\mu$, $\nu\leq\mu$.

\vspace{2 mm}

(iii) $\mathcal{L}(S_\nu)\subset\mathcal{L}(S_\mu)$, $\chi(S_\nu)\subset\chi(S_\mu)$ and therefore the characteristic foliations are also (inversely) nested hence the Cartan class of the successive derived systems tends to decrease, stopping only upon reaching integrability.
\vspace{2 mm}

\noindent
We can also assume, factoring if necessary $(P_\ell,S)$ modulo its characteristics, that

\vspace{2 mm}

(iv) The characteristic system $\chi(S)$ has maximum rank equal to $dim~P_\ell$. 

\vspace{2 mm}

\noindent
Factoring as well each $(P_\ell,S_\nu)$, modulo its own characteristics, we find a system $\overline{S}_\nu$ defined, as mentioned earlier, on a more suitable ambient space $P_{\ell-\nu}$ where it becomes, in analogy with \textit{S}, of maximum class (zero characteristics).
\vspace{2 mm}

\noindent
We shall now replace the fibration (5) by the \textit{structured} fibration

\begin{equation}
(P_\ell,S)\longrightarrow (P_{\ell-1},\overline{S}_1)\longrightarrow\cdots\longrightarrow (P_0,\overline{S}_\ell) ~,
\end{equation}

\vspace{1 mm}

\noindent
that now acquires, in its own right, a genuine structure of a \textit{prolongation space}\cite{Kumpera1999}. 
Let us now assume that \textit{S} is a flag system. In this case, we shall be surprised to observe that, 
besides the natural \textit{projection} functor, we gain much more namely, a \textit{prolongation} functor that ultimately provides accessibility to any term of a flag system directly from its Darboux foundation
since, as we shall see in the next section, the bottommost non-trivial derived system $S_{\ell-1}$ always 
identifies with the Darboux flag in $\bf{R^3}$. It should however be observed that, in the sequel, we  
simply try to provide the shortest and most geometrical approach to the \textit{merihedric} prolongations, when  restricted to the present context, and only transcribe what Élie Cartan already wrote in full generality more than 
a hundred years ago \cite{Cartan1914}.

\vspace{2 mm}

\newtheorem{lie}[TheoremCounter]{Theorem}
\begin{lie}
Every $\bf{i.a.}~\xi$ of $~\overline{S}_\mu$ prolongs (lifts) uniquely to an $\bf{i.a.}$ denoted $~p_{\mu-\nu}\xi$ of $~\overline{S}_\nu$ that projects onto $\xi$. Composition of two such prolongations is again a prolongation at the appropriate level and each prolongation map $~p_{\mu-\nu}:\mathcal{L}(\overline{S}_\mu)\longrightarrow\mathcal{L}(\overline{S}_\nu)$ is bijective, its inverse being the projection (extended to vector fields). Moreover, $p_\kappa$ is a Lie functor preserving all the operations including brackets and multiplication by functions and commutes with equivalences.
\end{lie}

\vspace{2 mm}

\noindent
As for the proof, it suffices of course to restrict our attention to \textit{S} and $S_1$, set $rank~S=\ell$, where $\ell$ is the \textit{length} of \textit{S}, consider the derived system written in the realm of its specific characteristic variables (or, equivalently, as being the inverse image of $\overline{S}_1$), and finally choose an adapted pseudo-normal form representation $\{\omega^1,~\cdots~,\omega^{\ell-1},\omega^\ell\}$ in such a way that $\{\omega^1,~\cdots~,\omega^{\ell-1}\}$ be in pseudo-normal form for the system $S_1$ and written with respect to a complete set of its characteristic variables $(x^1,~\cdots~,x^{\ell+1})$, the last form $\omega^\ell$ being the only one depending upon the last variable $x^{\ell+2}$. This being so, let $\xi=\sum A_i~\partial_i,~1\leq i\leq\ell+1$, be an infinitesimal automorphism of $~\overline{S}_1$ and let us determine an infinitesimal automorphism $\zeta=\sum A_i~\partial_i+f~\partial_{\ell+2}$ of \textit{S} that projects onto $\xi$, this last requirement forcing the above expression for $\zeta$. Consequently, the following congruences must hold:

\begin{equation}
\theta(\zeta)\omega^i\equiv 0~mod~\{\omega^1,~\cdots~,\omega^\ell\}, \hspace{3 mm} 1\leq i\leq\ell.
\end{equation}

\noindent
The first $\ell-1$ congruences obviously hold since $\xi$ is an infinitesimal automorphism of $~\overline{S}_1$ and $\theta(\partial_{\ell+2})\omega^i=\iota(\partial_{\ell+2})\omega^i=0$ for $i\leq\ell-1$. As for the last congruence involving $\theta(\zeta)\omega^\ell$, the pseudo-normal forms for \textit{S} enables us to write

\begin{equation*}
\theta(\zeta)\omega^\ell=\theta(\sum A_i~\partial_i)(dx^p+(x^{\ell+2}+c)~dx^q)~+
\end{equation*}
\begin{equation*}
+~\theta(f~\partial_{\ell+2})(dx^p+(x^{\ell+2}+c)~dx^q),~p,q\leq\ell+1,
\end{equation*}

\vspace{2 mm}

\noindent
hence

\begin{equation}
\theta(\zeta)\omega^\ell=dA_p+(x^{\ell+2}+c)~dA_q+f~dx^q
\end{equation}

\vspace{2 mm}

\noindent
and therefore, writing down explicitly the desired congruence in terms of all the $\omega^i$, the function \textit{f} becomes entirely \textit{and uniquely} determined by the $A_i$'s hence also by $\xi$. Another argument, in proving uniqueness, consists in showing, with the help of the above equation, that the only "vertical" vector field $f\partial_{\ell+2}$ i.e., tangent to the one dimensional characteristics of $S_1$ and that is an infinitesimal symmetry of \textit{S} necessarily vanishes. 

\vspace{2 mm}

It is interesting to rewrite the above expressions (and with the obvious notations) for the Engel and Cartan flags:

\begin{equation*}
p\xi=\xi+[(\partial C/\partial x^2~+~x^4\partial A/\partial x^2)x^3~+~(\partial C/\partial x^3~+~x^4\partial A/\partial x^3)x^4~+
\end{equation*}
\begin{equation*}
-~(\partial C/\partial x^1+~x^4\partial A/\partial x^1)]\partial_4\in\mathcal{L}(E),
\end{equation*}

\noindent
where $\xi\in\mathcal{L}(D)$ and $\xi=A~\partial_1+B~\partial_2+C~\partial_3$,

\begin{equation*}
\hspace{1 mm}p\xi=\xi+[(\partial D/\partial x^2~+~x^5\partial A/\partial x^2)x^3~+~(\partial D/\partial x^3~+~x^5\partial A/\partial x^3)x^4~+                         
\end{equation*}
\begin{equation*}
+~(\partial D/\partial x^4)x^5~-~(\partial D/\partial x^1~+~x^5\partial A/\partial x^1)]\partial_5\in\mathcal{L}(C),
\end{equation*}

\noindent
where $\xi\in\mathcal{L}(E)$ and $\xi=A~\partial_1+B~\partial_2+C~\partial_3+D~\partial_4$, and finally

\begin{equation*}
\hspace{6 mm}p\xi=\xi+[x^5(\eta-\partial D/\partial x^4)]\partial_5\in\mathcal{L}(CE),~where~\xi\in\mathcal{L}(E).\hspace{24 mm}
\end{equation*}

\vspace{4 mm}

\noindent
We thus see that the Darboux infinitesimal symmetries operate, via prolongation, as infinitesimal symmetries of any flag system and a simple inspection shows that this action is transitive for the Engel and the Cartan flags as well as for all the Cartan flags of arbitrary length $\ell$ (we recall that $\mathcal{L}(D)$ operates transitively on Darboux). However, it stops being so for all the remaining flags since, obviously, singularities must be preserved. We finally mention the following somewhat repetitive statement.

\vspace{2 mm}

\newtheorem{liefunctor}[CorollaryCounter]{Corollary}
\begin{liefunctor}
Let S be a flag system of length $\ell$. Each $\bf{i.a.}$ of S projects onto an $\bf{i.a.}$ of $\overline{S}_\nu,~1 \leq \nu \leq \ell$  and, conversely, every $\bf{i.a.}$ of $\overline{S}_\nu$, $\nu < \ell$, prolongs (lifts) to an $\bf{i.a.}$ of S. This prolongation-projection procedure commutes and establishes thereafter, for all $\nu,~0 < \nu < \ell$, a bijective Lie functor with source $\mathcal{L}(\overline{S}_\nu)$ and target $\mathcal{L}(S)$.
\end{liefunctor}

\vspace{2 mm}

Since we can freely shuffle around the derived systems, let us definitely forget about the upper bar on top of the $S_\nu$ and re-locate those systems at their most suitable places according to the needs of the context. Nevertheless, some care must be taken. For instance, the $\bf{i.a.}$ algebra of $S_\nu$, considered as a system on $P_\ell$, is the direct sum of the $\bf{i.a.}$ algebra of \textit{S} with the algebra of all characteristic automorphisms of $S_\nu$ \textit{i.e.}, the set of vector fields tangent to the characteristic foliation of $S_\nu$. 

\vspace{2 mm}

\noindent
A flag system of length $\ell$ and, in particular, a local model will be represented by a local basis $\{\omega^1,\cdots,\omega^\ell\}$ in such a way that $S_\nu$ is generated by $\{\omega^1,\cdots,\omega^{\ell-\nu}\}$ and is written with respect to its appropriate characteristic variables. In particular, $S_{\ell-1}$ is generated by a Darboux 1-form $\omega=\omega^1$ of maximum odd class equal to three and, for $\ell\geq 2$ , $S_{\ell-2}$ is generated by $\{\omega^1,\omega^2\}$, these forms being always representable by the (unique) Engel model. The previous lemma states (with the aforementioned precaution) that $\mathcal{L}(S)\simeq\mathcal{L}(S_{\ell-1})$ hence is also equivalent to $\mathcal{L}(D)$ (\textit{D=Darboux}).

\vspace{5 mm}

\noindent
\textbf{Infinitesimal isotropies of the Darboux flag}

\vspace{2 mm}

\noindent
We restrict our attention to the Darboux flag on $\bf{R^3}$ generated by the form $\omega=dx^2+x^3~dx^1$. Then

\begin{equation}
\mathcal{L}(D)=\{\xi\in\chi(\textbf{R}^3) : \theta(\xi)\omega\equiv 0~mod~\omega\}
\end{equation}

\vspace{2 mm}

\noindent
where $\theta(\xi)$ is the Lie derivative. There is a well known Hamiltonian characterization of $\mathcal{L}(D)$ that can be described as follows.

\noindent
Let $\mathcal{F}$ be the ring of $C^\infty$ functions on $\bf{R^3}$ and \textit{H} the contact Hamiltonian defined by

\begin{equation}
H:\xi\in\mathcal{L}(D)~\longrightarrow~\iota(\xi)\omega\in\mathcal{F}
\end{equation}

\vspace{2 mm}

\noindent
This map is bijective, the inverse being defined by

\begin{equation}
H^{-1}(f)=\xi=\partial f/\partial x^3~\partial/\partial x^1+(f-x^3~\partial f/\partial x^3)~\partial/\partial x^2 +
\end{equation}
\begin{equation*}
-(\partial f/\partial x^1-x^3~\partial f/\partial x^2)~\partial/\partial x^3
\end{equation*}

\vspace{2 mm}

\noindent
and, of course, it is determined by the Darboux flag \textit{S} up to a multiplicative factor depending upon the choice of the generating contact form $\omega$. We observe that $\theta(\xi)\omega=(\partial f/\partial x^2)\omega$ hence $\xi$ is an $\bf{i.a.}$ of the form $\omega$ if and only if $\partial f/\partial x^2=0$. This form being fixed, the Lie bracket of $\mathcal{L}(D)$ identifies, via \textit{H}, with the Lagrange bracket on $\mathcal{F}$ defined by $[f,g]=\xi g-g\zeta f$, where $\zeta=H^{-1}(1)$. Since the function \textit{f} in (11) is arbitrary, we infer the

\vspace{2 mm}

\newtheorem{lagrange}[LemmaCounter]{Lemma}
\begin{lagrange}
The algebra $\mathcal{L}(D)$ operates transitively on $\bf{R^3}$.
\end{lagrange}

\vspace{2 mm}

We next consider the $k-th$ order isotropy sub-algebra at the origin

\begin{equation}
\mathcal{L}_k(D)=\{\xi\in\mathcal{L}(D)~:~j^k_0~\xi=0\}
\end{equation}

\vspace{2 mm}

\noindent
where $j^k_0$ denotes the \textit{k-jet} at the origin. Then, on account of the expression (11), this sub-algebra identifies with

\begin{equation}
\mathcal{I}_k(D)=\{f\in\mathcal{F}~:~T_kA=T_kB=T_kC=0\}
\end{equation}

\vspace{2 mm}

\noindent
where

\begin{equation*}
H^{-1}(f)=A~\partial/\partial x^1+B~\partial/\partial x^2+C~\partial/\partial x^3,
\end{equation*}
\begin{equation*}
A=\partial f/\partial x^3,\hspace{3 mm}B=f-x^3~\partial f/\partial x^3,\hspace{3 mm}C=\partial f/\partial x^1-x^3~\partial f/\partial x^2
\end{equation*}

\vspace{2 mm}

\noindent
and where $T_kF$ denotes the $k-th$ order Taylor polynomial of \textit{F} at the origin. In particular,

\begin{equation*}
\mathcal{I}_0(D)=\{f\in\mathcal{F}~:~f(0)=f_1(0)=f_3(0)=0\}
\end{equation*}

\noindent
and
\begin{equation*}
\hspace{10 mm}\mathcal{I}_1(D)=\{f\in\mathcal{F}~:~f(0)=f_1(0)=f_2(0)=f_3(0)=f_{11}(0)=
\end{equation*}
\begin{equation*}
\hspace{23 mm}\hspace{10 mm}=f_{12}(0)=f_{13}(0)=f_{23}(0)=f_{33}(0)=0\}\hspace{2 mm}
\end{equation*}

\vspace{2 mm}

\noindent
where the lower indices denote partial derivatives of \textit{f}. Rather lengthy calculations will eventually show that

\begin{equation*}
\mathcal{I}_2(D)=\{f\in\mathcal{F}~:~ \cdots=f_{22}(0)=f_{111}(0)=f_{112}(0)=f_{113}(0)=
\end{equation*}
\begin{equation*}
\hspace{9 mm}=f_{122}(0)=f_{123}(0)=f_{133}(0)=f_{223}(0)=f_{233}(0)=f_{333}(0)=0\}
\end{equation*}

\vspace{2 mm}

\noindent
where "$\cdots$" denotes the equations defining $\mathcal{I}_1(D)$,

\begin{equation*}
\mathcal{I}_3(D)=\{f\in\mathcal{F}~:~ \cdots=f_{222}(0)=f_{1111}(0)=f_{1112}(0)=f_{1113}(0)=
\end{equation*}
\begin{equation*}
\hspace{4 mm}=f_{1122}(0)=f_{1123}(0)=f_{1133}(0)=f_{1222}(0)=f_{1223}(0)=f_{1233}(0)=\hspace{2 mm}
\end{equation*}
\begin{equation*}
\hspace{1 mm}=f_{1333}(0)=f_{2223}(0)=f_{2233}(0)=f_{2333}(0)=f_{3333}(0)=0\}\hspace{11 mm}
\end{equation*}

\vspace{2 mm}

\noindent
and
\begin{equation*}
\mathcal{I}_4(D)=\{f\in\mathcal{F}~:~ \cdots=f_{2222}(0)=f_{11111}(0)=f_{11112}(0)=f_{11113}(0)=
\end{equation*}
\begin{equation*}
=f_{11122}(0)=f_{11123}(0)=f_{11133}(0)=f_{11222}(0)=
\end{equation*}
\begin{equation*}
=f_{11223}(0)=f_{11233}(0)=f_{11333}(0)=f_{12222}(0)=
\end{equation*}
\begin{equation*}
=f_{12223}(0)=f_{12233}(0)=f_{12333}(0)=f_{13333}(0)=
\end{equation*}
\begin{equation*}
=f_{22223}(0)=f_{22233}(0)=f_{22333}(0)=f_{23333}(0)=f_{33333}(0)=0\}.
\end{equation*}

\vspace{3 mm}

\noindent
In general,

\begin{equation}
\mathcal{I}_k(D)=\{f\in\mathcal{F}~:~f_\alpha(0)~,~0\leq\mid\alpha\mid\leq k+1~,~\alpha\neq(0,k+1,0)\}
\end{equation}

\vspace{2 mm}

\noindent
where $\alpha=(\alpha_1,\alpha_2,\alpha_3)$ and $f_\alpha=\partial^{\mid\alpha\mid}f/(\partial x^1)^{\alpha_1}(\partial x^2)^{\alpha_2}(\partial x^3)^{\alpha_3}$.

\vspace{2 mm}

The number of independent equations defining $\mathcal{I}_k(D)$ is an invariant of the Darboux flag called the \textit{co-rank} of $\mathcal{I}_k(D)$ in $\mathcal{F}$~. In the case of $\mathcal{I}_0(D)$~, it is equal to the dimension of the orbit of the origin 0 in $\bf{R^3}$ under the infinitesimal action of $\mathcal{L}(D)$ and in the case of $\mathcal{I}_k(D)$~, it is equal to the dimension of an orbit under the infinitesimal action of the prolongation of $\mathcal{L}(D)$ to a suitable prolongation space of order \textit{k} e.g., the orbit of any Cartesian $k-th$ order frame, belonging to the fiber issued from the origin, under the infinitesimal action of the standard prolongation of $\mathcal{L}(D)$ to the 
$k-th$ order frame bundle. We observe, nevertheless, that frame bundles are not always the most suitable prolongation spaces to be looked for in the present context and, in most cases, the jet bundles and higher order Grassmannians are preferable. 

\vspace{4 mm}

\noindent
\textbf{Infinitesimal isotropies of the Engel flag}

\vspace{2 mm}

\noindent
We first show how the \textbf{i.a.} of the Darboux flag prolong uniquely onto the \textbf{i.a.} of the Engel flag. Let

\vspace{2 mm}

\begin{equation}
\xi=A~\partial_1+B~\partial_2+C~\partial_3+D~\partial_4\hspace{4 mm}(\partial_i=\partial/\partial x^i)
\end{equation}

\vspace{2 mm}
\noindent
be a vector field on $\bf{R^4}$ and let us assume that $\xi$ is an $\bf{i.a.}$ of the Engel flag \textit{E}. Then $\xi$ is also an $\bf{i.a.}$ of the derived system $E_1$ , an \textit{uplifted D (D=Darboux)}, whose characteristics are generated by the direction $\partial/\partial x^4$, hence the coefficients \textit{A, B} and \textit{C} of its projection in $\bf{R}^3$ only depend upon the variables ($x^1, x^2, x^3$). We therefore consider a vector field $\xi$ defined on $\bf{R^4}$ and assume that
\begin{equation*}
(i)~A,~B~and~C~only~depend~on~the~variables~(x^1,~x^2,~x^3)~and
\end{equation*}
\begin{equation*}
(ii)~the~projection~A~\partial_1+B~\partial_2+C~\partial_3~,~in~\bf{R^3},\hspace{25 mm}
\end{equation*}
\begin{equation*}
is~a~symmetry~of~the~Darboux~flag.\hspace{29 mm}
\end{equation*}

\vspace{2 mm}

\noindent
The condition $\theta(\xi)\omega^1\equiv 0~mod~\{\omega^1,~\omega^2\}$ is obviously fulfilled on account of (\textit{ii}) and the condition $\theta(\xi)\omega^2\equiv 0~mod~\{\omega^1,~\omega^2\}$ transcribes by

\vspace{2 mm}

\begin{equation*}
dC+x^4~dA+D~dx^1= \lambda(dx^2+x^3dx^1)~+~\mu(dx^3+x^4dx^1)
\end{equation*}

\vspace{2 mm}

\noindent
We thus obtain the equations
\begin{equation*}
\partial C/\partial x^1~+~x^4\partial A/\partial x^1~+~D=\lambda x^3~+~\mu x^4~,
\end{equation*}
\begin{equation*}
\partial C/\partial x^2~+~x^4\partial A/\partial x^2=\lambda~,\hspace{6 mm}
\end{equation*}
\begin{equation*}
\partial C/\partial x^3~+~x^4\partial A/\partial x^3=\mu~,\hspace{6 mm}
\end{equation*}

\vspace{2 mm}

\noindent
that provide
\begin{equation}
D=(\partial C/\partial x^2~+~x^4\partial A/\partial x^2)x^3~+~(\partial C/\partial x^3~+~x^4\partial A/\partial x^3)x^4~+
\end{equation}
\begin{equation*}
-~(\partial C/\partial x^1+~x^4\partial A/\partial x^1)
\end{equation*}

\vspace{2 mm}

\noindent
and, as already stated before, implies that an $\bf{i.a.}~\xi$ of the Engel flag is uniquely determined by its projection onto $\bf{R}^3$, whereupon the algebra
\begin{equation}
\mathcal{L}(E)=\{\xi\in\chi(\textbf{R}^4)~:~\theta(\xi)\omega^1\equiv 0,~\theta(\xi)\omega^2\equiv 0~mod~\{\omega^1,\omega^2\}\}
\end{equation}

\vspace{2 mm}

\noindent
is equal to the prolongation, to $\bf{R}^4$, of $\mathcal{L}(D)$, the prolongation algorithm being defined by the expression (16). Though the nature of this equation depends upon the choice of $\omega^1$ and $\omega^2$, the prolongation algorithm  itself is intrinsic. We can therefore identify $\mathcal{L}(E)$ with $\mathcal{L}(D)$ hence, ultimately, with $\mathcal{F}$ and, consequently, the $k-th$ order isotropy, at the origin, of the Engel flag becomes a proper sub-algebra of $\mathcal{I}_k(D)$. In other terms, $\mathcal{I}_k(E)\subset\mathcal{
I}_k(D)\subset\mathcal{F}$.

\vspace{2 mm}

\newtheorem{Engel}[LemmaCounter]{Lemma}
\begin{Engel}
The algebra $\mathcal{L}(E)$ operates transitively on $\bf{R^4}$.
\end{Engel}

\vspace{2 mm}

Further calculations show that
\begin{equation*}
\mathcal{I}_0(E)=\{f~:~f(0)=f_1(0)=f_3(0)=f_{11}(0)=0\},\hspace{16 mm}
\end{equation*}

\begin{equation*}
\mathcal{I}_1(E)=\{\cdots~=f_2(0)=f_{11}(0)=f_{12}(0)=f_{13}(0)=f_{23}(0)=
\end{equation*}
\begin{equation*}
\hspace{20 mm}=f_{33}(0)=f_{111}(0)=f_{112}(0)=f_{113}(0)=0\},
\end{equation*}

\begin{equation*}
\mathcal{I}_2(E)=\{\cdots~=f_{1111}(0)=f_{1112}(0)=f_{1113}(0)=f_{1122}(0)=
\end{equation*}
\begin{equation*}
=f_{1123}(0)=f_{1133}(0)=0\},\hspace{5 mm}
\end{equation*}

\begin{equation*}
\hspace{6 mm}\mathcal{I}_3(E)=\{\cdots~=f_{11111}(0)=f_{11112}(0)=f_{11113}(0)=f_{11122}(0)=
\end{equation*}
\begin{equation*}
\hspace{11 mm}=f_{11123}(0)=f_{11133}(0)=f_{11222}(0)=
\end{equation*}
\begin{equation*}
\hspace{17 mm}=f_{11223}(0)=f_{11233}(0)=f_{11333}(0)=0\},
\end{equation*}

\vspace{2 mm}

\noindent
and, in general,
\begin{equation*}
\mathcal{I}_k(E)=\{f~:~f_\alpha(0)=0~,~0\leq\mid\alpha\mid\leq k~+~1~,~\alpha\neq(0,k~+~1,0)~,
\end{equation*}
\begin{equation*}
\hspace{20 mm}~together~with~\alpha=(2,\alpha_2,\alpha_3)~,~\mid\alpha\mid=k~+~2\}.
\end{equation*}

\vspace{4 mm}

\noindent
\textbf{Infinitesimal isotropies of the Cartan homogeneous 3-flag}

\vspace{2 mm}

\noindent
We now show how the \textbf{i.a.} of the Engel flag, hence also of the Darboux flag, prolong uniquely onto the \textbf{i.a.} of the Cartan flag. Let
\begin{equation*}
\xi=A~\partial_1+B~\partial_2+C~\partial_3+D~\partial_4+E~\partial_5
\end{equation*}

\vspace{2 mm}

\noindent
denote a vector field on $\bf{R^5}$ and let us assume that
\begin{equation*}
A~\partial_1+B~\partial_2+C~\partial_3+D~\partial_4
\end{equation*}

\vspace{2 mm}

\noindent
is an \textbf{i.a.} of the Engel flag. Therefore $\xi=A~\partial_1+B~\partial_2+C~\partial_3$ is, in its own right, an \textbf{i.a.} of the Darboux flag hence the above assumption entails, by necessity, that $A=A(x^1, x^2, x^3), B=B(x^1, x^2, x^3), C=C(x^1, x^2, x^3)$ and $D=D(x^1, x^2, x^3, x^4)$. We now determine the coefficient \textit{E} in such a way that $\xi$ becomes an \textbf{i.a.} of the homogeneous Cartan model of length three. The condition
\begin{equation}
\theta(\xi)\omega^3\equiv 0~mod~\{\omega^1,\omega^2,\omega^3\}
\end{equation}

\vspace{2 mm}

\noindent
translates by the equation
\begin{equation*}
\theta(\xi) \omega^3=dD+x^5~dA+E~dx^1=\lambda \omega^1~+~\mu \omega^2~+~\eta \omega^3~,
\end{equation*}

\vspace{2 mm}

\noindent
hence by
\begin{equation*}
\partial D/\partial x^1~+~x^5\partial A/\partial x^1~+~E=\lambda x^3~+~\mu x^4~+~\eta x^5,
\end{equation*}
\begin{equation*}
\partial D/\partial x^2~+~x^5\partial A/\partial x^2=\lambda,\hspace{18 mm}
\end{equation*}
\begin{equation*}
\partial D/\partial x^3~+~x^5\partial A/\partial x^3=\mu,\hspace{18 mm}
\end{equation*}
\begin{equation*}
\hspace{4 mm}\partial D/\partial x^4=\eta.
\end{equation*}

\vspace{2 mm}

\noindent
The resulting expression for \textit{E} is therefore equal to
\begin{equation*}
E=\lambda x^3~+~\mu x^4~+~\eta x^5~-~(\partial D/\partial x^1~+~x^5 \partial A/\partial x^1)=
\end{equation*}
\begin{equation*}
=(\partial D/\partial x^2~+~x^5\partial A/\partial x^2)x^3~+~(\partial D/\partial x^3~+~x^5\partial A/\partial x^3)x^4~+
\end{equation*}
\begin{equation*}
+~(\partial D/\partial x^4) x^5~-~(\partial D/\partial x^1~+~x^5\partial A/\partial x^1),
\end{equation*}

\vspace{2 mm}

\noindent
and is determined by the coefficients \textit{A} and \textit{D} hence also by \textit{A,~B} and \textit{C}. In order to calculate the isotropy algebras of the Cartan flag, we assume that $T_kA=T_kB=T_kC=T_kD=0$ which means that the equations defining $\mathcal{I}_k(E)$ are satisfied, and then impose the condition $T_kE=0$. We find, eventually, the equations
\begin{equation*}
\mathcal{I}_0(C)=\{f~:~f(0)=f_1(0)=f_3(0)=f_{11}(0)=f_{111}(0)=0\},
\end{equation*}

\begin{equation*}
\mathcal{I}_1(C)=\{\cdots~=f_{1111}(0)=f_{1112}(0)=f_{1113}(0)=0\},\hspace{12 mm}
\end{equation*}

\vspace{2 mm}

\noindent
and
\begin{equation*}
\mathcal{I}_2(C)=\{\cdots~=f_{11111}(0)=f_{11112}(0)=f_{11113}(0)=\hspace{13 mm}
\end{equation*}
\begin{equation*}
\hspace{24 mm}=f_{11122}(0)=f_{11123}(0)=f_{11133}(0)=0\},\hspace{9 mm}
\end{equation*}

\vspace{2 mm}

\noindent
where the $"~\cdots~"$ denote the previous equations as well as those of the corresponding Engel isotropy. In general,
\begin{equation}
\mathcal{I}_k(C)=\{f~:~f_\alpha(0)=0~,~0\leq\mid\alpha\mid\leq k+1~,~\alpha\neq(0,k+1,0)~,
\end{equation}
\begin{equation*}
\mid\alpha\mid=k+2,~\alpha=(2,\alpha_2,\alpha_3)~,~\mid\alpha\mid=k+3~,~\alpha=(3,\alpha_2,\alpha_3)\}.
\end{equation*}

\vspace{4 mm}

\noindent
\textbf{Infinitesimal isotropies of the Cartan non-homogeneous 3-flag}

\vspace{2 mm}

\noindent
Let us show next how the \textbf{i.a.} of the Engel flag also prolong uniquely to the exceptional flag of length three. In this case, the condition on the coefficient \textit{E} writes
\begin{equation*}
\theta(\xi) \omega^3=dA+x^5~dD+E~dx^4=\lambda \omega^1~+~\mu \omega^2~+~\eta \omega^3,
\end{equation*}

\vspace{2 mm}

\noindent
hence
\begin{equation*}
x^5\partial D/\partial x^4~+~E=\eta x^5,\hspace{6 mm}
\end{equation*}
\begin{equation*}
\hspace{5 mm}\partial A/\partial x^1~+~x^5\partial D/\partial x^1=\lambda x^3~+~\mu x^4~+~\eta,
\end{equation*}
\begin{equation*}
\partial A/\partial x^2~+~x^5\partial D/\partial x^2=\lambda,\hspace{19 mm}
\end{equation*}
\begin{equation*}
\partial A/\partial x^3~+~x^5\partial D/\partial x^3=\mu,\hspace{19 mm}
\end{equation*}

\vspace{2 mm}

\noindent
and therefore, $E=x^5(\eta-\partial D/\partial x^4)$. Further calculations show that
\begin{equation}
\mathcal{I}_k(CE)=\mathcal{I}_k(E)
\end{equation}

\vspace{2 mm}

\noindent
for all integer \textit{k} , where \textit{CE} reads the initials of the Cartan non-homogeneous (Exceptional) flag and it is therefore quite surprising that, in the passage from Engel to Exceptional Cartan, the isotropy algebras stagnate. The  equality of the isotropies up to order 2 enabled Cartan to state his Criterion for the Monge equations in a non-homogeneous (and non-linear) context since, by a pure magic only accessible to the most illuminated, he managed to \textit{sneak} backstage Engel's homogeneity \cite{Cartan1910}. It is out of question to speak of dimensions when referring to isotropies since these are infinite dimensional algebras of Fréchet type.

\vspace{4 mm}

\noindent
\textbf{Infinitesimal isotropies of the Cartan homogeneous and nonhomogeneous $\ell$-flag}

\vspace{2 mm}

\noindent
Induction on the integer $\ell$ shows that the $k-th$ order isotropy algebra of the homogeneous Cartan flag of length $\ell$ is equal to (c.f., sect. 4)
\begin{equation*}
\mathcal{I}_k(C_\ell)=\{\mathcal{I}_k(C_{\ell-1}),~f_\alpha(0)=0~:~\mid\alpha\mid=k~+~\ell~,~\alpha=(\ell,\alpha_2,\alpha_3)\}~,
\end{equation*}

\vspace{2 mm}

\noindent
where, in the above expression, $\mathcal{I}_k(C_{\ell-1})$ represents all the equations resulting from $C_{\ell-1}$, hence
\begin{equation}
\mathcal{I}_k(C_\ell)=\{f~:~f_\alpha(0)=0~,~0\leq\mid\alpha\mid\leq k~+~1~,~\alpha\neq(0,k~+~1,0)~,
\end{equation}
\begin{equation*}
\alpha=(j,\alpha_2,\alpha_3)~,~\mid\alpha\mid=k~+~j~,~2\leq j\leq\ell\}.
\end{equation*}

\vspace{2 mm}

\noindent
We remark that $C_1=D$ and $C_2=E$.

\noindent
As for the exceptional (primary) flag of length $\ell~$, we find that
\begin{equation}
\mathcal{I}_k(CE_\ell)=\mathcal{I}_k(C_{\ell-1}).
\end{equation}

\vspace{4 mm}

We could keep on calculating isotropies for the other models made explicit in the previous section but prefer to have a closer look at the isotropy co-rank evolution resulting from successive prolongations, in various directions, of the algebra $\mathcal{L}(D)$ and thereafter the evolution of $\mathcal{L}(S)$ when the length of \textit{S} augments. We shall then conclude that 1 strictly diminishes isotropy i.e. there is an \textit{expanding} process in evolution, 3 maintains isotropy i.e., the process \textit{slows down} and finally $\bf{2}$ offers both alternatives. In the first case, it means that $\mathcal{L}(D)$ is progressively approaching, via homogeneity, a \textit{free} action whereas, in the second case, it maintains the same \textit{restraints} detectable in \textit{D}, \textit{E} and \textit{C}. Whenever 3 interferes in the evolution of 1, then simply the approach towards a free action is delayed. It is therefore very interesting to observe that this pattern somehow portrays the expansion of the universe where all mass and energy comes from the initial \textit{Big Bang}. And the \textit{black holes}, where are they? Perhaps the answer is given by some combination of \textbf{2}'s along the way as in the \textit{supermassive Sagittarius A}.

\vspace{2 mm}

Let \textit{S} be an arbitrary flag for which we know a $model~\{\omega^1~,~\cdots~,~\omega^\ell\}$ and where the first $\ell-\nu$ forms generate the $\nu$-\textit{th} derived system $S_\nu$. Then, in order to study the successive relations involving the isotropies, we take the last form $\varpi^2=\omega^{\ell-\nu}$ in $S_\nu$ as well as the last form $\varpi^1=\omega^{\ell-\nu+1}$ in $S_{\nu-1}$. Setting these two forms $\{\varpi^1,\varpi^2\}$ one on top of the other, as was done with the models, one can guess immediately, up to a re-indexation of the variables, whether the extension (or transition) of $S_\nu$ to $S_{\nu-1}$ is of the homogeneous type (Engel or Cartan type) or whether it is of the exceptional (singular) type (\textit{CE} type). There are, nevertheless, other types of singularities revealed by the non-vanishing, at the origin, of certain coefficients as exhibited by the fourth model of length-4 flags (see also \cite{Kumpera1982}.

\vspace{2 mm}

The general behaviour of flag systems is as follows. Whenever the extension of $S_\nu~$ to $~S_{\nu-1}$ is accomplished in the homogeneous Engel or Cartan way, the isotropy algebras strictly diminish (it is rather delicate to mention any sort of dimension) and, whenever this extension is singular of the type \textit{CE}, the isotropies stagnate (remain the same). However, when the singularity present in the transition is other than the two mentioned previously (e.g., the non-vanishing of some coefficient at the origin), both phenomena can actually take place and, moreover, such a singular transition can have a different nature from the \textit{CE} case. One can expect the unexpected from flag systems like, for example, the equivalence of the third model of a 4-flag, listed previously, with the standard write up of (3.1.) \cite[pg.227, fourth model]{Kumpera1982} but, instead of entering into details \textit{and} calculations, we find it more interesting and enlightening to illustrate this behavior by drawing a \textit{spider web}.
In order to do so, we first have to place \textit{tags} on the flags and, for that matter, we mention rapidly Mormul's tagging \cite{Mormul2000}. Since flags of rank one and two are always reducible to \textit{Darboux} and \textit{Engel}, it suffices to start tagging from rank three onward. For the homogeneous extension i.e., \textit{C}, we just write (1.) and for the exceptional one i.e., \textit{CE}, we write (3.). Finally, for the singular transition due to the presence of a non-vanishing coefficient at the origin, we write ($\bf{2}$.). With this in mind, \textit{C}=(1.), \textit{CE }=(3.), the first rank-4 model tags (1.1.), the second (1.3.), the third (3.1.), the fourth (3.$\bf{2}$.) and the fifth (3.3.). We can now deploy the first few knots of the web. A "segment" indicates that the two models have the same isotropy (i.e. stagnation) whereas a middle arrow signals the shrinking of the  isotropy towards the appointed target.

\begin{sidewaysfigure}[htp!]
    \center
    \includegraphics[width=.8\textwidth]{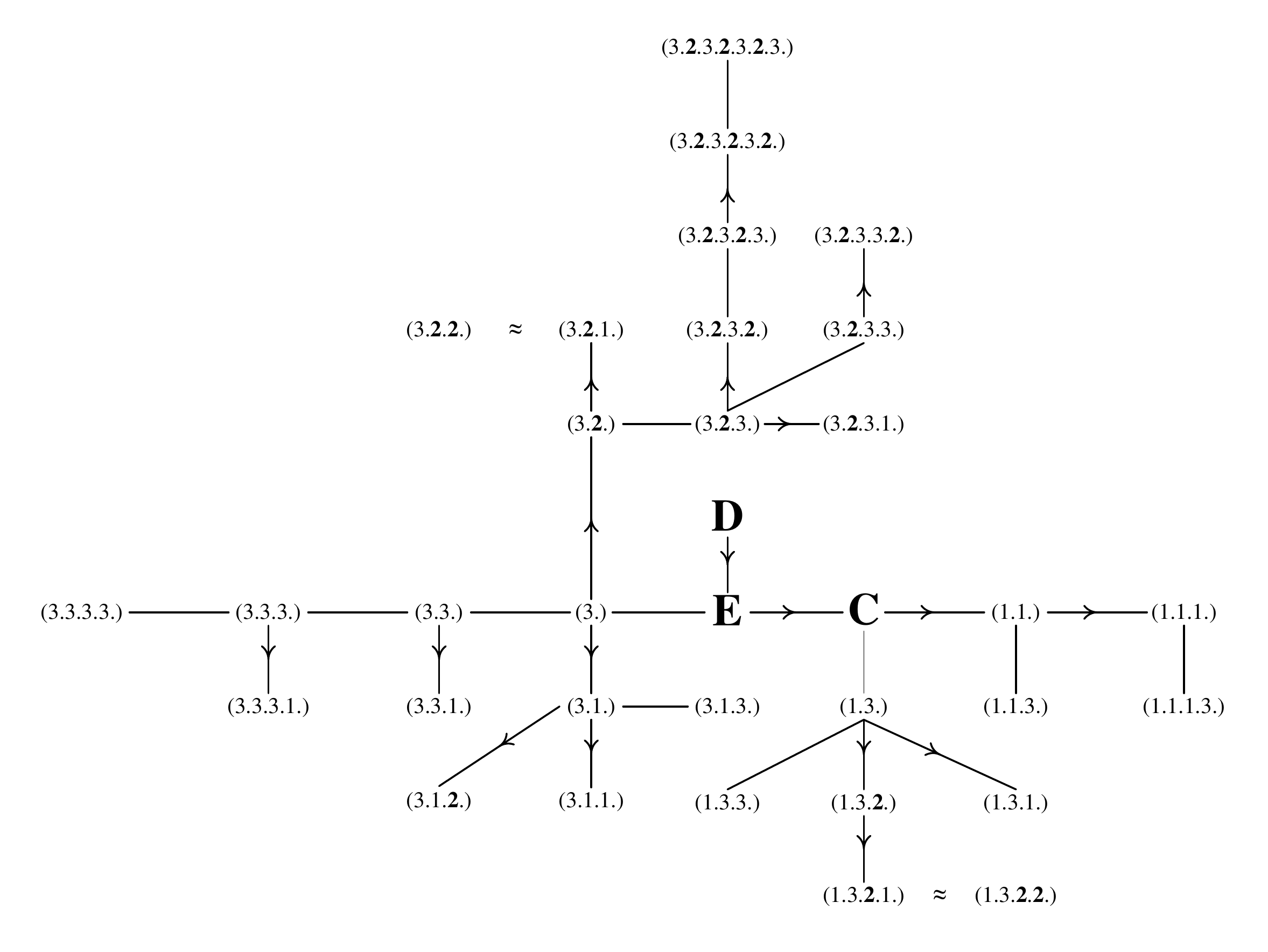}
    \caption{\textit{Spider Web}}
    \label{fig:diagram}
\end{sidewaysfigure}

\section{A non-linear approach}

Let \textit{S} be any Pfaffian system of rank $\ell$ defined on the manifold \textit{M}. We denote by $\mathcal{G}^k$ the groupoid of all $k-$jets of local automorphisms (local equivalences) of the system \textit{S} and by $\overline{\mathcal{G}}^k$ the groupoid of all $k-$jets of local diffeomorphisms that induce $(k-1)-st$ order jet equivalences of \textit{S}. When \textit{k} = 1, $\overline{\mathcal{G}}^1$ is the set of all the $1-$jets $X=j^1_x\varphi$ such that $X^*(S_y)=\varphi^*(S_y)=S_x$, where $x=\alpha(X)$ is the source of \textit{X} and $y=\beta(X)$ is the target. When \textit{k}~>~1, $\overline{\mathcal{G}}^k$ is the set of all the \textit{k}-jets $X=j^k_x\varphi$ such that $j^{k-1}_x(\varphi^*(S))=j^{k-1}_x(S)$, where \textit{S} is considered as a section of the Grassmannian of all the $\ell$-planes in $T^*M$. Since the construction of the characteristic system $\chi(S)$ associated to \textit{S} requires one differentiation, the elements of $\overline{\mathcal{G}}^k$ only determine $(k-2)-d$ order jet equivalences of $\chi(S)$. Furthermore, since the germ of a flag system of length $\ell$ is $(\ell+1)-$\textit{determined}, we can state the following

\vspace{2 mm}

\newtheorem{Handel}[LemmaCounter]{Lemma}
\begin{Handel}
Let S be a flag system of length $\ell$. Then $\mathcal{G}^k=\overline{\mathcal{G}}^k,~ k>\ell$. 
\end{Handel}

\vspace{2 mm}

This lemma provides a finite algorithm for computing $\mathcal{G}^k$ when $k>\ell$. According to the Lemma 2, any model of length $\ell$ defined on $\bf{R}^{\ell+2}$ is locally equivalent, in the neighborhood of the origin, to the restriction of an \textit{elementary} model to a neighborhood of a suitable point $p\in\bf{R}^{\ell+2}$. Consequently, in order to study the equivalence classes of germs of flag systems of length $\ell$, it suffices to study the orbits of the groupoids $\overline{\mathcal{G}}^{\ell+1}$ associated to the elementary models of that same length. However, it should be mentioned that the orbits of a given elementary model can represent other non-equivalent elementary models. For the transitive model there is just one orbit (more precisely, all the elements in the orbit are equivalent) hence it can only represent, locally, the transitive model. For the exceptional model in dimension 5, there are just two orbits hence it can represent, locally, either the transitive model or the exceptional one. Let us now take, for example, the model (3.3.) in dimension 6. A simple calculation (re-scaling of the coordinates) shows that the orbits of the corresponding groupoid $\overline{\mathcal{G}}^{4+1}$ represent, locally, the models (1.1.), (1.3.), (3.$\bf{2}$.) but cannot represent (3.1.) as can also be checked independently by verifying that, whatever the point in $\bf{R}^6$, either the \textit{small growth vector} or the \textit{reduced tensors} as well as the \textit{reduced nilpotent algebras} cannot be equal to those of (3.1.) at the origin. As for the \textit{good behavior} of the above groupoids, the answer is negative. Of course, the groupoids corresponding to the transitive models of any length are as good as possible namely, they are regularly embedded sub-manifolds of the corresponding jet spaces and are integrable, formally integrable, involutive (considered as PDE) and each higher order groupoid is the prolongation of the preceding one \cite{Kumpera1972}. As for the others, they all must preserve the singularities present in each model. For example, the $\beta$-projection of the $\alpha$-fiber at the origin is, in the case of the model (3.), its singular locus namely, the hyper-plane $x^5=0$. A similar and sometimes even more complex situation arises for the other non-transitive models and their groupoids will not comply to Ehresmann's requirements \cite{Ehresmann1951,Ehresmann1958}. For the sake of catching at least a twinkle of light at the end of the tunnel, we shall, in the next section, write down the equations defining the groupoids $\overline{\mathcal{G}}^1$ for the first few flags.

\section{Equations defining $\overline{\mathcal{G}}^1$ for the initial flags}
Calculations are, for many, a boring activity though, for others, become indispensable. Let us therefore exhibit the equations defining the first order equivalence groupoids $\overline{\mathcal{G}}^1$ for Darboux, Engel, Cartan, Cartan Exceptional and a few other flags. Since these equations just involve first order derivatives, they only take into account the \textit{fibre to fiber} equivalence of \textit{S} and ignore the associated characteristic system. We first write down the full set of equations for \textit{D} and thereafter indicate
only those equations that have to be added due to the transition to a longer model. Moreover, since the calculations involved are excessively long though very simple and straightforward, we just indicate the results. Let us denote by
$(y^j)$ and $(x^i)$ the dependent and independent variables respectively i.e., $\varphi:(x^i)\longmapsto(y^j)$.

\vspace{4 mm}

\noindent
$\bf{The~Darboux~flag.}$

\vspace{2 mm}

\noindent
We have to find the equations that translate the following congruence:
\begin{equation*}
\varphi^*\omega=dy^2+y^3~dy^1=\alpha(dx^2+x^3~dx^1)~,~\omega=dx^2+x^3~dx^1.
\end{equation*}

\vspace{2 mm}

\noindent
and obtain
\begin{equation}
y^2_3=-y^3y^1_3,\hspace{8 mm}y^2_1=-y^3y^1_1+x^3(y^2_2+y^3y^1_2).
\end{equation}

\vspace{2 mm}

\noindent
$\bf{The~Engel~flag.}$

\vspace{2 mm}

\noindent
Here we have to transcribe, by means of equations, the congruence (transition equations)
\begin{equation*}
\varphi^*\omega^2=dy^3+y^4~dy^1=\alpha(dx^2+x^3~dx^1)+\beta(dx^3+x^4~dx^1)
\end{equation*}

\vspace{2 mm}

\noindent
and obtain, besides the equations (23),
\begin{equation}
y^3_1=-y^4y^1_1+x^3(y^3_2+y^4y^1_2)+x^4(y^3_3+y^4y^1_3).
\end{equation}

\vspace{2 mm}

\noindent
$\bf{The~Cartan~homogeneos~(1.)~flag.}$

\vspace{2 mm}

\noindent
The Cartan homogeneous (transitive) context adds the following (transition) congruence
\begin{equation*}
\varphi^*\omega^3=dy^4+y^5~dy^1=\alpha(dx^2+x^3~dx^1)+\beta(dx^3+x^4~dx^1)+\gamma(dx^4+x^5~dx^1)
\end{equation*}

\vspace{2 mm}

\noindent
where after we obtain the equation
\begin{equation}
y^4_1=-y^5y^1_1+x^3(y^4_2+y^5y^1_2)+x^4(y^4_3+y^5y^1_3)+x^5y^4_4.
\end{equation}

\vspace{2 mm}

\noindent
$\bf{The~Cartan~homogeneous~(1.1.)~flag.}$
\begin{equation*}
y^5_1=-y^6y^1_1+x^3(y^5_2+y^6y^1_2)+x^4(y^5_3+y^6y^1_3)+x^5y^5_4+x^6y^5_5.
\end{equation*}

\vspace{2 mm}

\noindent
$\bf{The~Cartan~homogeneous~(1.1.1.)~flag.}$
\begin{equation*}
y^6_1=-y^7y^1_1+x^3(y^6_2+y^7y^1_2)+x^4(y^6_3+y^7y^1_3)+x^5y^6_4+x^6y^6_5+x^7y^6_6.
\end{equation*}

\vspace{2 mm}

\noindent
$\bf{The~Cartan~homogeneos~\ell-flag~(1.1.~\cdots~.1.).}$
\begin{equation*}
y^{\ell+1}_1=-y^{\ell+2}y^1_1+x^3(y^{\ell+1}_2+y^{\ell+2}y^1_2)+x^4(y^{\ell+1}_3+y^{\ell+2}y^1_3)+x^5y^{\ell+1}_4+
\end{equation*}
\begin{equation*}
x^6y^{\ell+1}_5+x^7y^{\ell+1}_6+~\cdots~+x^{k+1}y^{\ell+1}_k+~\cdots~+x^{\ell+1}y^{\ell+1}_\ell+x^{\ell+2}y^{\ell+1}_{\ell+1}.
\end{equation*}

\vspace{2 mm}

\noindent
$\bf{The~Cartan~exceptional~(3.)~flag.}$

\vspace{2 mm}

\noindent
In this case we start with the (transition) congruence
\begin{equation*}
\varphi^*\omega^3=dy^1+y^5~dy^4=\alpha(dx^2+x^3~dx^1)+\beta(dx^3+x^4~dx^1)+\gamma(dx^1+x^5~dx^4)
\end{equation*}

\vspace{2 mm}

\noindent
and find the equation
\begin{equation}
y^5y^4_4=x^5[y^1_1+y^5y^4_1-x^3(y^1_2+y^5y^4_2)-x^4(y^1_3+y^5y^4_3)].
\end{equation}

\vspace{2 mm}

\noindent
Since $y^4_4\neq 0$ ($\varphi$ being an automorphism, it projects "downwards" and, accordingly, its components depend on "lesser" variables - see sect. 5), we infer that, on the singular hyperplane $\{x^5=0\}$, also hold the equalities $y^5=y^5_j=0,~j\neq 5,$ and consequently that $(3.)\not\simeq (1.)$ in any neighborhood of the origin.

\vspace{2 mm}

\noindent
$\bf{The~Cartan~exceptional~(3.3.)~flag.}$
\begin{equation}
y^6y^5_5=x^6[y^4_4+y^6y^5_4-x^5(y^4_1+y^6y^5_1-x^3(y^4_2+y^6y^5_2)+
\end{equation}
\begin{equation*}
-x^4(y^4_3+y^6y^5_3))].
\end{equation*}

\vspace{2 mm}

\noindent
$\bf{The~flag~(3.1.).}$
\begin{equation*}
\omega^4=dx^5+x^6~dx^4,\hspace{4 mm}\varphi^*\omega^4=d\varphi^5+\varphi^6~d\varphi^4,
\end{equation*}

\vspace{2 mm}

\noindent
leads to the equation
\begin{equation}
y^6y^4_4=-y^5_4+x^6y^5_5+x^5[y^5_1+y^6y^4_1-x^3(y^5_2+y^6y^4_2)-x^4(y^5_3+y^6y^4_3)],
\end{equation}

\vspace{2 mm}

\noindent
that reduces, on the singular locus $\{x^5=x^6=0\}$, to $y^5=y^5_j=y^6=y^6_j=0,~j~<~5$. We infer, once again, that $(3.1.)\not\simeq(3.\bf{2})$.

\vspace{2 mm}

\noindent
$\bf{The~flag~(3.1.3.).}$
\begin{equation*}
\omega^5=dx^4+x^7~dx^6,\hspace{4 mm}\varphi^*\omega^5=d\varphi^4+\varphi^7~d\varphi^6,
\end{equation*}

\vspace{2 mm}

\noindent
and leads to the equation
\begin{equation*}
y^7y^6_6=x^7[y^4_4+y^7y^6_4-x^6y^7y^6_5-x^5\{y^4_1+y^7y^6_1-x^3(y^4_2+y^7y^6_2)+
\end{equation*}
\begin{equation*}
-x^4(y^4_3+y^7y^6_3)\}].
\end{equation*}

\vspace{2 mm}

\noindent
Since $y^6_6\neq 0$, we infer that, on the singular locus $\{x^7=0\}$, also hold the equalities $y^7=y^7_j=0,~j\neq 7$.

\vspace{2 mm}

\noindent
$\bf{The~flag~(3.1.3.1.).}$
\begin{equation*}
\omega^6=dx^7+x^8~dx^6,\hspace{4 mm}\varphi^*\omega^6=d\varphi^7+\varphi^8~d\varphi^6,
\end{equation*}

\vspace{2 mm}

\noindent
and leads to the equation
\begin{equation*}
y^8y^6_6=-y^7_6+x^8y^7_7+x^7\epsilon=
\end{equation*}
\begin{equation*}
=-y^7_6+x^8y^7_7+x^7[y^7_4+y^8y^6_4-x^6(y^7_5+y^8y^6_5)+
\end{equation*}
\begin{equation*}
-x^5\{y^7_1+y^8y^6_1-x^3(y^7_2+y^8y^6_2)-x^4(y^7_3+y^8y^6_3)\}].
\end{equation*}

\vspace{2 mm}

\noindent
where the coefficient $\epsilon$ is determined by the equation (one of seven equations)
\begin{equation*}
y^7_4+y^8y^6_4=x^5\gamma+x^6\delta+\epsilon
\end{equation*}

\vspace{2 mm}

\noindent
resulting from the requirements imposed by the congruence equation (transition equation).

\vspace{2 mm}

\noindent
$\bf{The~flag~(3.1.3.3.).}$
\begin{equation*}
\omega^6=dx^6+x^8~dx^7,\hspace{4 mm}\varphi^*\omega^6=d\varphi^6+\varphi^8~d\varphi^7,
\end{equation*}

\vspace{2 mm}

\noindent
leads to the equation
\begin{equation*}
y^8y^7_7=x^8[y^6_6+y^8y^7_6-x^7\epsilon]=
\end{equation*}
\begin{equation*}
=x^8[y^6_6+y^8y^7_6-x^7\{y^6_4+y^8y^7_4-x^5(y^6_1+y^8y^7_1+
\end{equation*}
\begin{equation*}
-x^3(y^6_2+y^8y^7_2)-x^4(y^6_3+y^8y^7_3))-x^6(y^6_5+y^8y^7_5)\}],
\end{equation*}

\vspace{2 mm}

\noindent
the $\epsilon$ coefficient resulting from the equations

\vspace{2 mm}

$\hspace{17 mm}y^6_4+y^8y^7_4=x^5\gamma+x^6\delta+\epsilon$ and $y^6_6+y^8y^7_6=x^7\epsilon+\phi$.

\vspace{2 mm}

\noindent
On the singular locus $\{x^8=0\}$, the first equations reduce to $y^8=y^8_j=0,~j\neq 8$. 

\vspace{2 mm}

\noindent
$\bf{The~flag~(3.1.3.1.3.).}$
\begin{equation*}
\omega^7=dx^6+x^9~dx^8,\hspace{4 mm}\varphi^*\omega^7=d\varphi^6+\varphi^9~d\varphi^8,
\end{equation*}

\vspace{2 mm}

\noindent
leads to the equation
\begin{equation*}
y^9y^8_8=x^9[y^6_6+y^9y^8_6-x^8(y^9y^8_7)-x^7\epsilon]=
\end{equation*}
\begin{equation*}
=x^9[y^6_6+y^9y^8_6-x^8(y^9y^8_7)-x^7\{y^6_4+y^9y^8_4-x^6(y^6_5+y^9y^8_5)+
\end{equation*}
\begin{equation*}
-x^5(y^6_1+y^9y^8_1-x^3(y^6_2+y^9y^8_2)-x^4(y^6_3+y^9y^8_3))\}].
\end{equation*}

\vspace{2 mm}

\noindent
On the singular locus $\{x^9=0\}$, these equations reduce to $y^9=y^9_j=0,~j\neq 9$, and $y^9_9y^8_8=\psi$, the coefficients $\epsilon$ and $\psi$ being determined by the equations

\vspace{2 mm}

$\hspace{3 mm}y^6_4+y^9y^8_4=x^5\gamma+x^6\delta+\epsilon$, $y^6_6+y^9y^8_6=x^7\epsilon+x^8\phi+\psi$ and $y^9y^8_8=x^9\psi$.

\vspace{2 mm}

\noindent
$\bf{The~flag~(3.1.3.3.1.).}$
\begin{equation*}
\omega^7=dx^8+x^9~dx^7,\hspace{4 mm}\varphi^*\omega^7=d\varphi^8+\varphi^9~d\varphi^7,
\end{equation*}

\vspace{2 mm}

\noindent
leads to the equations
\begin{equation*}
y^9y^7_7=-y^8_7+x^9y^8_8+x^8\phi=
\end{equation*}
\begin{equation}
=-y^8_7+x^9y^8_8+x^8[y^8_6+y^9y^7_6-x^7\{y^8_4+y^9y^7_4+
\end{equation}
\begin{equation*}
-x^5(y^8_1+y^9y^7_1-x^3(y^8_2+y^9y^7_2)-x^4(y^8_3+y^9y^7_3))-x^6(y^8_5+y^9y^7_5)\}],
\end{equation*}

\vspace{2 mm}

\noindent
$\phi$ being determined by $y^8_6+y^9y^7_6=x^7\epsilon+\phi$ and $y^8_7+y^9y^7_7=x^8\phi+x^9y^8_8$.

\vspace{3 mm}

Just to exemplify the usefulness of the above calculations, let us prove for example that
\begin{equation*}
(3.\textbf{2}.3.3.\textbf{2}.)~\not\simeq~(3.\textbf{2}.3.3.\textbf{2-}.)
\end{equation*}
where the minus sign indicates that the constant in this pseudo-model is negative.
\nocite{Cheaito1998}\nocite{Hilbert1912}\nocite{Kumpera1974}\nocite{Kumpera1991}\nocite{Kumpera1985}\nocite{Mormul1999}

Considering the last elementary model (3.1.3.3.1.), it amounts to show that the groupoid $\mathcal{G}^k$ associated to this
model does not contain any element with
\begin{equation*}
source~x=(0,0,0,0,0,1,0,0,1)~and~target~y=(0,0,0,0,0,1,0,0,-1).
\end{equation*}
The equation (29) reduces at the point \textit{x} to $y^9y^7_7=-y^8_7+x^9y^8_8$ and, by (28), $y^8_7=0$ hence, for $x^9=1$ and $y^9=-1,~y^7_7=-y^8_8$ at the point \textit{x}. Differentiating, with respect to $x^8$, the last equation resulting from the congruence requirement in the (3.1.3.3.) model, namely the equation $y^8y^7_7=x^8\phi$, we obtain $y^7_7y^8_8=\phi+x^8\partial\phi/\partial x^8$ hence, at the point $x,~y^7_7y^8_8=\phi~<~0$. However, by the second equation in (27), $y^6_6=\phi~<~0$ and, for $x^5=0$, the expression (26) reduces to $y^6y^4_4=x^6y^5_5$ hence to $y^4_4=y^5_5\neq 0$ at the point $x~(x^6=y^6=1)$. Differentiating the above expression with respect to $x^6$, we obtain $y^6_6y^4_4=y^5_5~(y^4$ and $y^5$ do not depend on $x^6)$ and consequently $y^6_6=1~>~0$.

Inasmuch and using a similar argument, we can also prove the non-equivalences (c* and $\overline{c}^{*}$ indicating the constants at the corresponding \textbf{2} positions):
\begin{equation*}
(3.\textbf{2}.3.\textbf{2}.3.\textbf{2}(c^{10}).)~\not\simeq~(3.\textbf{2}.3.\textbf{2}.3.\textbf{2}(\overline{c}^{
10}).)\hspace{5 mm} whenever~c^{10}\neq\overline{c}^{10},
\end{equation*}
\begin{equation*}
(1.3.1.1.\textbf{2}.\textbf{2}.\textbf{2}(c^{11}).)~\not\simeq~(1.3.1.1.\textbf{2}.\textbf{2}.\textbf{2}(\overline{c}^{11}).)\hspace{5 mm} whenever~c^{11}\neq\overline{c}^{11}.
\end{equation*}

\vspace{2 mm}

\noindent
and conclude, in particular, that there exist \textit{uncountably} many non-equivalent flag systems. Moreover, in case non-vanishing constants do appear in the models with a specific given code, their non-equivalence classes can eventually be parametrized by these constants, which is the case above. Such an approach can be found in \cite{Mormul2005} though treated in a very different setting and only in a rather specific case where the 
number 3 does not appear in the code \textit{i.e.}, the flag being therefore free of any inversion.

\vspace{2 mm}

\noindent
However, we can do much better in trying to apply \textit{Lie Theory} to the equivalence groupoids arising in connection with flag systems and believe that much benefit can result. A serious drawback however was, for a long time, the fact that Lie Theory required strong regularity conditions that unfortunately are not fulfilled in the present context. The answer came with Malgrange \cite{Malgrange1970,Malgrange1972} who showed that one can in fact always argue \textit{up to a closed Zariski subset} and this fits perfectly well in our setting. In fact, the groupoids here considered are everywhere very nice sub-manifolds of jet spaces except above the \textit{singular loci} of the flags. However, these loci are either hyperplanes, or else, linear subspaces of lower dimension of numerical spaces hence obviously closed and Zariski. Since
the elements of these groupoids are jets of local equivalences, we infer that if the source of an element belongs to the singular locus then, forcefully, the target will also belong to it and a similar statement can be repeated for the target. Denoting by $\Lambda$ the above locus, the following relations hold: For any $x,y\in\Lambda,~\beta(\alpha^{-1}(x))=\alpha
(\beta^{-1}(y))=\Lambda$ and therefore the singular locus in the groupoid is nothing else but the union of the inverse images $\alpha^{-1}(x)$ or $\beta^{-1}(y)$ of any two given points as well as the inverse image of $\Lambda\times\Lambda$ by $\alpha\times\beta$ and consequently is also closed and Zariski. We should perhaps undergo a more careful analysis since the singular locus can admit a stratification but, at present, this seems irrelevant. We also remark, at this point, that the equations given in previous examples for the non-homogeneous flags and defining these groupoids are \textit{singular i.e. not of constant rank} except of course in the (homogeneous) Cartan realm that repeats itself as well, for singular flags, in the complement of the above Zariski set. At present, we shall not elongate on the theme of \textit{Lie Equations} since much has been already written on this subject. Our purpose is exclusively practical and aims in showing how Lie theory can be used within the context of flag systems.

\vspace{2 mm}

As remarked previously, a flag system of length $\ell$ is $(\ell+1)-$determined and, this being so, we can calculate $\overline{\mathcal{G}}^{k+1},~k>\ell$ or its equivalent $\mathcal{G}^{k+1}$ by simply lifting, to order $k+1$, the equations defining the $k-th$ order groupoid, $k>\ell$, and add to them afterwards all the \textit{total derivatives} of the previous equations thus obtaining a set of equations defining the groupoid at order $k+1$. However, if we perform the same operation at levels $h\leq\ell$ we still obtain an $(h+1)-st$ order groupoid (e.g., $\overline{\mathcal{G}}^h$ and $\overline{\mathcal{G}}^{h+1}$) that however might be larger than the desired equivalence groupoid at that order and contain "non-equivalence defining" jets. This in fact, can be observed by calculating first prolongations of $\mathcal{G}^h$ \textit{via} total derivations and remarking that an additional condition of the type $x^3=0$ is still missing so as to obtain $\mathcal{G}^{h+1}$ in the case of Darboux and Engel. The calculations involved are of course unduly long, though revealing.

\section{On classification}

The intersections of the derived systems with the characteristic systems provide very useful information and it is worthwhile to recall the following result \cite{Montgomery2001} based on the inclusions $S_k\subset\chi(S_{k+1}),~0\leq k\leq\ell-2$ and where $S=S_0$ and $S_\ell=\chi(S_\ell)=0$.

\vspace{2 mm}

\newtheorem{Zhito}{Proposition}
\begin{Zhito}
The transitive Cartan model is the unique model of length $\ell$ for which $S_k\neq\chi(S_{k+2})$ at the origin, $0\leq k\leq\ell-3$. Furthermore, each sub-sequence for which the above inequality holds, equality holding for the remaining terms, characterizes all the elementary models of lenght $\ell$.
\end{Zhito}

\vspace{2 mm}

Curiously enough, we can also characterize all the elementary models of a given length in quite a different manner, but having much verbal resemblance with the above statement. As was already mentioned, each time that a transition proceeds according to the Cartan homogeneous fashion, the isotropy algebras strictly diminish (their functional co-rank strictly augments) whereas, if the transition is singular of \textit{CE} type, the isotropy remains the same. We thus infer that the nature of the $S_\nu$ isotropies, $\nu\leq\ell$, in what concerns their progressive "decay", provides a full classification of the elementary flags of length $\ell$ and further exhibits the corresponding models.
As for the non-elementary flags, there is some hope of a similar classification if we investigate continuous (or differentiable) deformations of such isotropy algebras.

\vspace{2 mm}

It is also worthwhile to observe that the intersections of derived systems with characteristic systems at different levels can provide further information about a given flag. These new Pfaffian systems are intrinsically associated to \textit{S} and can be either integrable or not. For example, the system $S\cap\chi(S_1)$ has rank equal to 4, is integrable for the models (3.1.) and (3.3.) but is non-integrable for (3.\textbf{2}.).\nocite{Cartan1893}

\section*{Acknowledgements}
Our acknowledgements are especially due to Rui Almeida and Jean Pradines for precious help and encouragement.  Inasmuch, we would also like to write an apology similar to that written by Édouard Goursat in his book with the only difference that, whereas he states \textit{presque à chaque page de ces trois chapitres}, we should claim \textit{presque à chaque ligne}.

\bibliographystyle{spmpsci}      
\bibliography{references}   

\end{document}